\documentclass{amsart}

\usepackage{amssymb}
\usepackage{amsmath}
\usepackage{amsthm}
\usepackage{amscd}
\usepackage{calc}
\usepackage{tikz}
\usepackage{pifont}
\usepackage{color}
\usepackage{xypic}
\usepackage{ulem}
\usetikzlibrary{cd}

\usepackage{graphicx}

\def\qed{$\s$}
\def\eps{\varepsilon}
\def\Gproj{\operatorname{G-proj}}
\def\proj{\operatorname{proj}}
\def\fin{\operatorname{fin}}
\def\deg{\operatorname{deg}}

\def\uGproj{\operatorname{\underline{G-proj}}}
\def\usub{\operatorname{\underline{sub}}}

\def\uH{\underline H}

\def\rank{\operatorname{rank}}
\def\Ab{{\mathcal A}b}
\def\mod{\operatorname{mod}}

\def\sub{\operatorname{sub}}
\def\Ker{\operatorname{Ker}}
\def\Im{\operatorname{Im}}
\def\Hom{\operatorname{Hom}}

\def\Ima{\operatorname{Im}}  

\def\lmatrix#1{(\begin{smallmatrix}0&0\\#1&0\end{smallmatrix})}
\renewcommand\bar{\overline}
\newcommand\boldit[1]{\textbf{\textit{#1}}}

\numberwithin{equation}{subsection}
\newtheorem{thm}[equation]{Theorem} 
\newtheorem{lem}[equation]{Lemma}

\newtheorem{prop}[equation]{Proposition}

\theoremstyle{definition}
\newtheorem{rem}[equation]{Remark}
\newtheorem{defin}[equation]{Definition}

\definecolor{darkgreen}{rgb}{0,0.5,0}

\newcommand{\s}{\hfill \blacksquare}

\newenvironment{smallpmatrix}{\left(\begin{smallmatrix}}{\end{smallmatrix}\right)}
\def\smallmat#1#2#3#4{\left(\begin{smallmatrix} #1\;&#2\\#3\;&#4\end{smallmatrix}\right)}

\input prepictex   \input pictex    \input postpictex
\begin{document}
\thispagestyle{empty}
\phantom m
\vspace{-1cm}

       {\tiny{\tt \jobname,} \today}
       \bigskip

\begin{center}
  {\large\bf  Gorenstein-Projective Modules over the Ring of Dual Integers }

  \medskip
  Xiu-Hua Luo
  and Markus Schmidmeier

\bigskip \parbox{10cm}
{\footnotesize{\bf Abstract:}
 
  The ring of dual integers is the bounded polynomial ring $\mathbb Z[\eps]=\mathbb Z[T]/(T^2)$
  with integer coefficients.  We describe the (finitely generated)
  Gorenstein-projective $\mathbb Z[\eps]$-modules
  as the torsionless $\mathbb Z[\eps]$-modules, while 
  the stable category of $\Gproj\mathbb Z[\eps]$ modulo projectives is shown to be 
  equivalent to the orbit category $\mathcal D^b(\mathbb Z)/[1]$ of the derived category
  of the integers.
  It follows that the latter carries the structure of a triangulated category.

  \smallskip
  The category $\Gproj\mathbb Z[\eps]$ is related to the embeddings of a subgroup
  in a free abelian group and has  a quotient which is equivalent to the category of
  finite abelian groups.  In fact, we present a cube which has as vertices eight related
  categories and as edges in each of the three directions functors which are
  related to push-down functors modulo the shift; canonical functors to stable categories;
  and homology functors, respectively. We note that
  in $\Gproj\mathbb Z[\eps]$ uniqueness of direct sum decomposition fails.
 
}

\medskip \parbox{10cm}{\footnotesize{\bf MSC 2010:}
  18G25 (relative homological algebra),
  20K27 (subgroups of abelian groups)
}

\medskip \parbox{10cm}{\footnotesize{\bf Key words:}
   Gorenstein-projective modules, dual integers, orbit category,
    push-down functor, cohomology functor,
    derived categories, torsionless modules, subgroup embeddings,
    Krull-Remak-Schmidt property. }

\end{center}

\section{Introduction}
\label{sec-intro}
By $\mathbb Z[\eps]$ we denote the ring of dual integers which is the bounded
polynomial ring $\mathbb Z[T]/(T^2)$ with integer coefficients.
In this paper, all modules are finitely generated right modules. 
Note that the classification of all modules (or even of all finite modules)
over $\mathbb Z[\eps]$
up to isomorphy is not considered a feasible problem, see Remark~\ref{rem-wild}. Hence our interest in this paper is in the Gorenstein-projective modules for which we
  discuss classification, properties and related categories. 

\subsection{The Gorenstein-projective modules}
\label{sec-one-one}
We focus on the category $\Gproj\mathbb Z[\eps]$ of
Gorenstein-projective modules over the ring
$\mathbb Z[\eps]$ of dual integers. 
Gorenstein-projective modules can be dated back to the
G-dimension zero modules which were introduced  by Auslander in \cite{Aus67}.
Then, their properties were intensively studied by many reseachers,
such as in \cite{AB69, EJ95, Buch, BGS87, CPST08}.
One of the fundamental problems is to explicitly construct
all the Gorenstein-projective modules over a given ring (see \cite{LiZ, Z11}).

\medskip
The category $\Gproj\mathbb Z[\eps]$ in this paper is understood
since it is equivalent to other categories, namely the
torsionless $\mathbb Z[\eps]$-modules
and the orbits of bounded
complexes of free abelian groups, see Theorem~\ref{four-cats-equ}. 
But we note that $\Gproj \mathbb Z[\eps]$
does not have the Krull-Remak-Schmidt property (decomposition into indecomposable objects
in general is not unique).

\subsection{The stable category of Gorenstein-projective modules}
\label{sec-one-two}
Of particular interest is the stable category
$\uGproj \mathbb Z[\eps]$ of Gorenstein-projective
$\mathbb Z[\eps]$-modules modulo projectives.

Here is our main result which is shown in Section~\ref{sec-stable}.

\begin{thm}
  \label{thm-main}
  The following categories are equivalent:
  \begin{enumerate}
  \item[(a)] The stable category of $\uGproj \mathbb Z[\eps]$ of Gorenstein-projective
    $\mathbb Z[\eps]$-modules modulo projectives.
  \item[(b)] The stable category $\usub(\mathbb Z[\eps]) $
    of all submodules of a finitely generated free $\mathbb Z[\eps]$-module   modulo projectives.     
  \item[(c)] The orbit category $\mathcal D^b(\mathbb Z)/[1]$
    of the derived category of the integers.
  
  \end{enumerate}
\end{thm}

In this paper, we focus on properties of Gorenstein-projective modules over the
ring of dual integers.

\medskip
 
  We note that the above theorem can also be established when we replace the ring of integers
  by any principal ideal domain, or more general, by any
   hereditary Noetherian integral domain $R$.
  This includes the Dedekind domains which occur as rings of algebraic integers in number theory  (see Remark~\ref{rem-number-rings}).

\medskip
Many mathematicians have studied the links
between $\uGproj\mathbb Z[\eps]$ and other important categories.
R.\ O.\ Buchweitz (\cite[Theorem 4.4.1]{Buch})
has shown that if $\Lambda$ is Gorenstein, then there is a triangle-equivalence
$D^b_{sg}(\Lambda)\cong\uGproj(\Lambda)$ (see 
also \cite[Theorem 4.6]{Hap2}). Recently, T. Honma and S. Usui \cite{HU} proved
for an arbitrary monomial algebra $\Lambda$ that
the stable category of graded Gorenstein-projective $\Lambda$-modules
$\uGproj^{\mathbb Z}(\Lambda)$ is triangle
equivalent to the bounded derived category of a path algebra of Dynkin type $\mathbb A$.

\medskip
The essential point of our theorem can be regarded as the
gradability of the Gorenstein-projective $R[\eps]$-modules.
Consider the
$\mathbb Z$-graded ring $R[\eps]$ whose degree zero part is $R$ and $\deg \eps=1$.
Regarding the action of $\eps$ as differential, then
$C^b(\proj R)$ coincides with the category $\Gproj^{\mathbb Z} R[\eps]$ of $\mathbb Z$-graded
Gorenstein-projective $R[\eps]$-modules by Theorem \ref{four-cats-equ},
and $K^b(\proj R)=D^b(\mod R)$ coincides the stable category
$\underline{\Gproj}^{\mathbb Z} R[\eps]$.
Since the degree shift functor $(1)$ of $\underline{\Gproj}^{\mathbb Z} R[\eps]$
coincides with the suspension $[1]$, one obtains an equivalence
$D^b(\mod R)/[1] \simeq (\underline{\Gproj}^{\mathbb Z} R[\eps])/(1)$, where the latter
category can be regarded as a full subcategory of the ungraded stable
category $\underline{\Gproj} R[\eps]$.

  \subsection{A triangulated structure for the orbit category of $\mathcal D^b(\mathbb Z)$}
  We recall from \cite[Theorem~1]{Ke} that for a finite dimensional hereditary algebra $A$,
  the category $\mathcal D^b(A)/[1]$ carries a natural triangulated structure.
  We obtain as a consequence of Theorem~\ref{thm-main}
  that also $\mathcal D^b(\mathbb Z)/[1]$ is a triangulated category, as follows.
  We note in Section~\ref{sec-Gproj} that$\Gproj\mathbb Z[\eps]$ is a Frobenius category.
  Hence the stable category $\uGproj\mathbb Z[\eps]$ is triangulated.  The categorical equivalence
  shown in Theorem~\ref{thm-main} yields that $\mathcal D^b(\mathbb Z)/[1]$ is a triangulated
  category.

  \smallskip
  It remains open if $\mathcal D^b(\mathbb Z)/[1]$ carries a triangulated structure
  as defined in \cite{Ke}
  and if, in this case, the equivalence $\uGproj \mathbb Z[\eps]\cong\mathcal D^b(\mathbb Z)/[1]$
  is a triangle functor.

\subsection{A chain of categories }
\label{sec-one-three}
It turns out that the two categories in Sections~\ref{sec-one-one} and \ref{sec-one-two}
occur as part of a chain of categories
linked by dense functors:  Let $\proj(\mathbb Z)$ be the category of finitely generated
projective $\mathbb Z$-modules, that is, of finitely generated free abelian groups.
The first category in the chain is the category $C^b=C^b(\proj(\mathbb Z))$
of bounded complexes with coefficients in $\proj(\mathbb Z)$;
the factor modulo the shift is the
category $\Gproj\mathbb Z[\eps]$ of Gorenstein-projective modules
considered in Section~\ref{sec-one-one};
we obtain as factor modulo the projective objects
the stable category $\uGproj\mathbb Z[\eps]$ in Section~\ref{sec-one-two}
which is equivalent to $\mathcal D^b(\mathbb Z)/[1]$;
finally we arrive at the category $\Ab_{\fin}$
of finite abelian groups by applying the
functor $H^\bullet_{\mathcal P,1}$ induced from the homology functor
$H^\bullet:C^b(\proj(\mathbb Z))\to \amalg \mod\mathbb Z$. 
$$C^b(\proj(\mathbb Z))\;\stackrel S\longrightarrow\;
\Gproj\mathbb Z[\eps]\;\stackrel{P_1}\longrightarrow\;
\mathcal D^b(\mathbb Z)/[1]\;\stackrel{H^{\bullet}_{\mathcal P,1}}\longrightarrow\;
\Ab_{\rm fin}$$

\subsection{A cube whose faces are commutative}
\label{sec-one-four}
As indicated in the above chain, we are dealing with three types of dense functors:
 Starting at the category $C^b=C^b(\proj(\mathbb Z))$,
  they are the push-down functor $S:C^b\to C^b/[1]\cong\Gproj\mathbb Z[\eps]$ modulo the shift onto
  the orbit category;
  the canonical functor $P: C^b\to C^b/\mathcal P\cong K^b(\mathbb Z)$
  modulo the ideal of all maps which factor through a projective object;
  and the functor $H^\bullet: C^b\to H^\bullet C^b\cong\amalg\mod\mathbb Z$ obtained by taking cohomology.
  Different functors of type $S$, $P$ and $H^\bullet$ can be composed and the product
  does not depend on the order. 
Thus, we obtain a cube with faces commutative squares (see Section~\ref{sec-cube}).

$$
\begin{tikzcd}
  3 & &  \Ab_{\fin} & \\
  & & & \\
  2 & \mathcal D^b(\mathbb Z)/[1] \arrow[uur,"H^{\bullet}_{ \mathcal P,1}"]
  & \amalg \Ab_{\fin} \arrow[uu,"S_{H, \mathcal P}",near start]
  \arrow[from=ddr,"P_H"',near start]
  & \mod\mathbb Z \arrow[uul,"P_{H,1}"'] \\
  & & &\\
  1 & K^b(\mathbb Z) \arrow[uu,"S_{\mathcal P}"]
  \arrow[uur,"H^\bullet_{\mathcal P}", near start]
  & \Gproj\mathbb Z[\eps] \arrow[uul,"P_1"',near start, crossing over]
  \arrow[uur, "H^\bullet_1", near start, crossing over]
  & \amalg \mod \mathbb{Z}
  \arrow[uu,"S_H"'] \\
  & & & \\
  0 & & C^b(\proj(\mathbb Z)) \arrow[uul,"P"] \arrow[uu,"S"]
  \arrow[uur,"H^{\bullet}"'] & 
\end{tikzcd}
$$

\medskip 
The path sketched in Section~\ref{sec-one-three} leads from $C^b(\proj(\mathbb Z))$ at the
bottom of the cube via $S$ to the
category $\Gproj\mathbb Z[\eps]$ in the center, then 
then via $P_1$ to the upper left, and finally via $H^\bullet_{\mathcal P,1}$
to $\Ab_{\fin}$ at the top.

\medskip

\section{The structure of Gorenstein-projective $\mathbb Z[\eps]$-modules}
\label{sec-2}
In this section, we study the structure of Gorenstein-projective modules over the
ring $\mathbb Z[\eps] = \mathbb Z[T ]/(T^2 )$ of dual integers.
We  show that the three categories in Theorem~\ref{four-cats-equ} are equivalent.
Some of our results hold more generally. 
In Section~\ref{sec-torsionless}, we describe the torsionless
$R[\eps]$-modules.
Then in Section~\ref{sec-Gproj}, we describe the Gorenstein-projective modules over the
ring $ \mathbb Z[\eps] $ as the torsionless $\mathbb Z[\eps]$-modules.
Finally, the relation with the orbit category for bounded complexes of finitely generated
free abelian groups is shown in Section~\ref{sec-two-four}.

\subsection{Torsionless $R[\eps]$-modules}
\label{sec-torsionless}

Suppose that $R$ is a hereditary Noetherian integral domain.
Each $R[\eps]$-module is a pair $(N, e)$ where $N\in\mod R$ and the $R$-morphism
$e: N\to N$ satisfies $e^2=0$.
If $N$ is projective, then
$(N, e)\cong (A\oplus B,\lmatrix \mu)$ where $\mu:A\to B$ a monomorphism
between projective $R$-modules (see Lemma \ref{isom}).
Thus, the object $(N,e)$ embeds into $(B\oplus B,\lmatrix1)$
and hence is a subobject of $R[\eps]^n$ for some natural number $n$.
Following \cite{RZ}, we call a submodule of a module of the form $R[\eps]^n$
\boldit{torsionless} and denote by $\sub(R[\eps])$
the full subcategory of $\mod R[\eps]$
of all torsionless $R[\eps]$-modules.

\medskip
The following lemma is an important observation which we use to show
the main result of this section.
 
\begin{lem}
  \label{isom}
  Let $R$ be a hereditary Noetherian integral domain.
  Each $R[\eps]$-module $(N, e)$ with $N$ being projective
  is isomorphic to
  $(I\oplus \Ker e, \smallmat 00\mu0)$,
    where $I$ is a subobject of $N$ isomorphic to $\Ima e$
    and $\mu: I\to \Ker e$ is the restriction of the multiplication by $e$.
\end{lem}

\begin{proof}
  Since $N$ is a projective object in $\mod R$ and
  the map $e:N\to N$ satisfies $e^2=0$,
  the image of $e$  can be embedded into the kernel of $e$,
  say $\Ima e\to \Ker e$.
  Since $R$ is hereditary, then $\Ima e$ and $\Ker e$ are projective.
  It follows that the short exact sequence
  $$ 0\to \Ker e\stackrel{\sigma}\longrightarrow N
  \stackrel{\pi}\longrightarrow \Ima e \to 0$$
  splits where $\sigma$ is the inclusion and $\pi$ the multiplication by $e$.
  So there exist morphisms $\sigma': N\to \Ker e$ and $\pi':\Ima e\to N$
  such that $\sigma'\sigma=1_{\Ker e}$, $\pi\pi'=1_{\Ima e}$ and
  $N=I\oplus\Ker e$ where $I=\Ima \pi'$.
  It is easy to check that the following square is commutative
  $$\xymatrix{N \ar[r]^-{e}
    \ar[d]_-{\left( \begin{smallmatrix} \pi'\pi \\ \sigma'\end{smallmatrix}\right)}
    & N \ar[d]^-{\left( \begin{smallmatrix} \pi'\pi \\ \sigma'\end{smallmatrix}\right)} \\
    I\oplus\Ker e
    \ar[r]^-{\lmatrix\mu}
    & I\oplus\Ker e}
  $$

  where $\mu:I\to \Ker e$ is the restriction of $e$.
  The map $\left( \begin{smallmatrix}\pi'\pi\\ \sigma'\end{smallmatrix}\right)$ is injective,
    hence an isomorphism. That is to say,
    $\left( \begin{smallmatrix}\pi'\pi\\ \sigma'\end{smallmatrix}\right):\
      (N, e)\to (I\oplus \Ker e,
      \left( \begin{smallmatrix} 0& 0 \\ \mu&0\end{smallmatrix}\right))$
        is an isomorphism of $R[\eps]$-modules.
\end{proof}

We conclude this subsection with a remark about homomorphisms.

\begin{rem}
  \label{rem-homomorphisms}
  Let $(A\oplus B,\lmatrix u)$
    and $(A'\oplus B',\lmatrix{u'})$
      be $R[\eps]$-modules.
   It is straightforward to check that 
    the homomorphisms are of the form 
    $\smallmat a0cb$ 
    satisfying
    $\smallmat a0cb
      \lmatrix u=\lmatrix{u'}
      \smallmat a0cb$.
      That is to say, there is a diagram 
      \[\xymatrix{A\ar[r]^-{u}\ar[d]_-
        {a}\ar[dr]^{c} & B \ar[d]^-{b} \\
        A'\ar[r]^-{u'} &B'}\]
      where the square commutes (although the two triangles may not commute).
\end{rem}

\subsection{The Gorenstein-projective $\mathbb Z[\eps]$-modules}
\label{sec-Gproj}
  \begin{defin}\label{def-Gp}
    Let $R$ be any ring and $\mod R$ be the category of all finitely generated
    right $R$-modules.
    An $R$-module $M$ is
    \boldit{Gorenstein-projective} if there is a double infinite exact
    sequence
    $$P^\bullet: \qquad \cdots\longrightarrow P^{-1}\longrightarrow P^0
    \stackrel{d^0}\longrightarrow P^1\longrightarrow P^2\longrightarrow\cdots$$
    of finitely generated projective $R$-modules such that the dual sequence $\Hom_R(P^\bullet,R)$
    is exact and
    $M\cong\Ker d^0$.
    A Gorenstein-projective module $M$ is
    \boldit{strongly Goren\-stein-projective} if all the modules and all the
    maps in the above sequence can be chosen to be equal.
  \end{defin}

  We denote by $\Gproj R$ the full subcategory of $\mod R$ consisting of all
  Gorenstein-projective modules.
    For $R$ a Noetherian ring, the category of finitely generated $R$-modules
    is abelian. We recall from \cite[Proposition~3.8]{Beli} that in this case,
    the category $\Gproj R$ of Gorenstein-projective $R$-modules is a
    Frobenius category.

\begin{prop}\label{prop-subR}
  Let $\mathbb Z[\eps]$ be the ring of dual integers.
  \begin{enumerate}
  \item   The category $\Gproj \mathbb Z[\eps]$
    of  Gorenstein-projective $\mathbb Z[\eps]$-modules
    is exactly the category $\sub(\mathbb Z[\eps])$. 
  \item Every Gorenstein-projective $\mathbb Z[\eps]$-module
    is strongly Gorenstein-projective.
  \end{enumerate}
\end{prop}

\begin{proof}
  (1) First, we need to show that the Gorenstein-projective $\mathbb Z[\eps]$-modules are exactly the objects in $\sub(\mathbb Z[\eps])$.
  In fact, $\Gproj\mathbb Z[\eps]$ is contained in
  $\sub(\mathbb Z[\eps])$ since the kernel of $d^0$ in the definition
  of Gorenstein-projective is a subobject of a finitely generated
  projective $\mathbb Z[\eps]$-module, hence of a free $\mathbb Z[\eps]$-module of finite rank.

  \smallskip
  Conversely, every object
  $(A\oplus B,\lmatrix u)$
    in $\sub(\mathbb Z[\eps])$ where $B$ is a free abelian group of finite rank and $u: A\to B$ an inclusion map,
    is Gorenstein-projective.
  Consider the short exact sequence
  \begin{eqnarray*}
    {{\mathcal E :}}&  0\longrightarrow
  (A\oplus B,\lmatrix u)
     & \stackrel{%
      \begin{smallpmatrix} -\eta\\ \iota\end{smallpmatrix}}
    \longrightarrow
    (A\oplus A,\lmatrix 1)
        \oplus (B\oplus B,\lmatrix 1)\\
          & & \stackrel{(\iota,\eta)}\longrightarrow
          (A\oplus B,\lmatrix u)\longrightarrow 0
  \end{eqnarray*}

  where the maps are given by
  $(a,b)\mapsto ((0,-a),(a,b))$ and $((r,s),(t,u))\mapsto (r,s+t)$.
  The sequence $P^\bullet$ in the definition is given by
  $$P^i=(A\oplus A,\lmatrix1)
    \oplus(B\oplus B,\lmatrix 1)
      \quad\text{and} \quad
  d^i=\begin{smallpmatrix} -\eta & 0\\ \iota & \eta\end{smallpmatrix},$$
  for every $i\in\mathbb Z$, where $d^i$ is the composition
  $\begin{smallpmatrix} -\eta\\ \iota\end{smallpmatrix}\circ (\iota,\eta)$
  of the two labeled maps in $\mathcal E$. It is straightforward to verify
  that the dual sequence $\Hom(P^\bullet,\mathbb Z[\eps])$ is also exact. 
  Hence, $(A\oplus B,\lmatrix u)$
    is Gorenstein-projective.

  Both $\Gproj\mathbb Z[\eps]$ and $\sub(\mathbb Z[\eps])$ are full subcategories of mod$\mathbb Z[\eps]$ with the same objects.
  Hence these two categories are not only equivalent but actually the same. 
  
  \medskip
  (2). From the exact sequence ${\mathcal E}$, we can compose a double infinite exact sequence $P^{\bullet}$ with the same projective modules and maps. Then by Definition~\ref{def-Gp}, this assertion is proved.
\end{proof}
 
\begin{rem} By the above Proposition \ref{prop-subR} and Remark \ref{rem-homomorphisms}, the Gorenstein-projective $\mathbb Z[\eps]$-modules can be described as the embeddings $A\stackrel{u}\to B$ between free abelian groups. Note that homomorphisms in $\Gproj \mathbb Z[\eps]$ have the form $\left(\begin{smallmatrix} a&0\\c&b \end{smallmatrix}\right)$,while
morphisms in the embedding category are given by commutative squares, hence are
just pairs of the form $(a, b)$. 

More generally, comparing the relations between Gorenstein-projective modules and monomorphism categories (see \cite{LiZ, LZ1, LZ2}), for the quiver  Q of Dynkin type $A_2$ and a 
$k$-algebra $\Lambda$, the Gorenstein-projective $\Lambda\otimes_k kA_2$-modules can be described as monomorphism categories $\mathcal S(\Gproj\Lambda)$: They have the form $A\stackrel{f}\to B$ where $A$ and $B$ are Gorenstein-projective $\Lambda$-modules and $f$ is a monomorphism in the category $\Gproj\Lambda$. 

\medskip
Let $\proj(\mathbb Z)$ be the category of finitely generated free abelian groups,
  that is, projective $\mathbb Z$-modules.
  The relation between $\Gproj \mathbb Z[\eps]$ and $\mathcal S(\proj(\mathbb Z))$
  will be given in Section~\ref{sec-three-two}.
\end{rem}

\subsection{The equivalence between $  \sub(R[\eps])$ and $ C^b(\proj(R))/[1]$ }
\label{sec-two-four}

Let $R$ be a hereditary Noetherian integral domain
and $\proj(R)$ the subcategory consisting of the projective objects in $\mod R$.
By $ C^b(\proj(R))$ we denote the full subcategory of $C^b( \mod R )$
  whose objects are the bounded complexes with modules in $\proj(R)$.
  When there is no ambiguity regarding the underlying ring $R$,
    we abbreviate $C^b=C^b(\proj(R))$. 
The following property of indecomposable complexes in $C^b$
is a key observation to reveal the relation between $ \sub(R[\eps])$ 
and $ C^b/[1]$.

\begin{prop}
  \label{prop-cb} Every object in  $C^b(\proj(R))$
  is a direct sum of objects of width at most 2
  of the form $\cdots\to 0\to U\to F\to 0\to \cdots $ where  $F\in \proj(R)$
  and $U\to F$ is the inclusion of a submodule.
\end{prop}

\begin{proof}
    Let $X^\bullet:\cdots \to 0\to X\stackrel d\to X'\stackrel{d'}\to X''
  \stackrel{d''}\to\cdots$ be a complex in $C^b$.
  Since $R$ is a hereditary ring and $X'$ is projective,
  the map $d:X\to\Im d$ is a split epimorphism.  Hence $X\cong \Ker d
  \oplus\Im  d$. Similarly, $X'\cong\Ker d'\oplus\Im d'$. Etc.
  We have the following commutative diagram:
  $$
  \begin{tikzcd} 0 \arrow[r] \arrow[d] & X \arrow[r,"d"]
    \arrow[d,"\cong"] & X' \arrow[r,"d'"] \arrow[d,"\cong"] &
    X'' \arrow[r] \arrow[d,"\cong"] & \cdots \\
    0 \arrow[r] & \Ker d\oplus\Im d
    \arrow[r,"\left({0\atop 0}{\iota\atop 0}\right)"] 
    & \Ker d'\oplus\Im d' \arrow[r,"\left({0\atop0}{\iota\atop 0}\right)"]
    & \Ker d''\oplus\Im d''\arrow[r] & \cdots
  \end{tikzcd}
  $$
  This shows that $X^\bullet$ is a direct sum of complexes of width two
  of the form $\Im d\stackrel\iota\to\Ker d'$,
  $\Im d'\stackrel\iota\to\Ker d''$, etc.
  As a direct summand of an object in $\proj(R)$, $\Ker d'$ is projective, hence free
  \cite[Theorem 14.6]{f1}.
\end{proof}
\medskip
We usually denote an orbit in $C^b/[1]$ by the
  complex which represents it, or, if it has width at most two, by
  a map or as an inclusion.

  \medskip
  For two objects $X,Y\in C^b$, the homomorphisms between the orbits
  are given by
  $$\Hom_{C^b/[1]}(X,Y)=\bigoplus_{i\in\mathbb Z}\Hom_{C^b}(X,Y[i]).$$

  Thus, for indecomposable objects
  $(A\stackrel{u}\hookrightarrow B)$, $(A'\stackrel{u}\hookrightarrow B')$,
  homomorphisms are given by three homomorphisms in
  $\mathcal C$, $a:A\to A'$, $b: B\to B'$
  and $c:A\to B'$, subject to the condition that $a,b$ make the square commutative.
  $$
  \begin{tikzcd}
    A \arrow[r, "{u}"] \arrow[d,"a"] & B\arrow[d,"b"]\\
    A' \arrow[r, "{u'}"] & B'
  \end{tikzcd}
  \qquad
  \begin{tikzcd}
    & A \arrow[r,"{u}"] \arrow[d,"c"] & B\\
    A' \arrow[r,"{u'}"] & B' &
  \end{tikzcd}
  $$
  Composition is given by the composition of morphisms of complexes.

  We define a functor $ \eta:\ C^b\to \sub(R[\eps])$
  sending a bounded complex
  $X: \ \cdots0\to X^i\stackrel{d^i_X }\longrightarrow X^{i+1}\to
  \cdots \stackrel{d^{j-1}_X }\longrightarrow X^j\to 0\cdots$
  to the $R$-module
$$\left(\bigoplus\limits_{t=i }\limits^{j}
  X^t, \left( \begin{smallmatrix} 0&0&\cdots&0&0\\ d^i_X&0&\cdots&0&0\\ & & \ddots&\ddots&  \vdots\\ 0&0&\cdots&d^{j-1}_X&0\end{smallmatrix}\right)\right).$$

    From Lemma \ref{isom}, we know that each $R$-module $(N, e)$ occurs in the image,
    $ \eta(\cdots 0\to I \stackrel{\mu }\longrightarrow \Ker e\to 0\cdots)
    \cong (N, e) $.
    Hence $ \eta$ is dense.
    It is obvious that for each $X\in C^b$ and $i\in \mathbb{Z}$,
    $\eta(X[i])= \eta (X)$. So $ \eta$ induces a dense functor
     $ \eta_1:\ C^b/[1]\to \sub(R[\eps])$.

\begin{lem}\label{orbit-diff}
  The induced functor $ \eta_1:\ C^b (\proj(R))/[1]\to \sub(R[\eps])$
  is a categorical equivalence.
\end{lem}

\begin{proof}
  By Proposition \ref{prop-cb}, each indecomposable object in
  $\ C^b$ has width at most 2. 
  Hence we can suppose without loss of generality that 
  \begin{eqnarray*} X: &
    \cdots \to 0\to X^1\stackrel{d^{1}_X }\longrightarrow X^2\to 0\to \cdots,\\
    Y: &
    \cdots \to 0\to Y^1\stackrel{d^{1}_Y }\longrightarrow Y^2\to 0\to \cdots
  \end{eqnarray*}
  are indecomposable with $  X^1, Y^1$ at the $1^{\rm st}$ branch and $d^{1}_X, d^{1}_Y$ being injective
  (for the stalk complex, we can take the $1^{\rm st}$ branch zero and the $2^{\rm nd}$
  branch an indecomposable projective object in $\mod R$.)
  \medskip
  
 We see that the homomorphism groups $\Hom_{C^b/[1]}(X,Y)$
  determined above and $\Hom_{\sub(R[\eps])}(\eta_1 X, \eta_1 Y)$ determined in
  Remark~\ref{rem-homomorphisms} are isomorphic.
\end{proof}
 
\begin{rem}
   Note that the push-down functor $\eta:C^b(\mod \mathbb Z)\to   \mod(\mathbb{Z}[\eps])  $
  on the full module category is not dense.  For example, the object
  $$(N,e)=(\mathbb Z/(p^2),\mu_p)$$
    where $\mu_p$ denotes the multiplication by $p$ does not occur in the image of $\eta$, since $N$ is indecomposable
    (but the object $(\mathbb Z/(p)\oplus\mathbb Z/(p),\lmatrix1)$ does occur in the image of $\eta$).

    The remark shows that free abelian groups are very special.
\end{rem}

From Propositon \ref{prop-subR} and Lemma \ref{orbit-diff},
we obtain one of the main results of this paper as follows:

\begin{thm}\label{four-cats-equ}
  Let  $\mathbb Z[\eps]$ be the ring of dual integers
  and $\proj(\mathbb Z)$ the category of finitely generated free abelian groups.
  The following categories are equivalent. 
  $$C^b(\proj(\mathbb Z))/[1]\; \cong\; \sub (\mathbb Z[\eps])
  \;\cong\; \Gproj\mathbb Z[\eps]$$  \qed
\end{thm}  

\section{Some properties of Gorenstein-projective $\mathbb Z[\eps]$-modules }
\label{sec-properties}
 
  In this section we show first that there is an epivalence
  from the category of Gorenstein-projective $\mathbb Z[\eps]$-modules
  to the category $\mathcal S(\proj\mathbb Z)$ of embeddings
  of a subgroup in a
  free abelian group (Proposition~\ref{prop-epivalence}).
  Embeddings provide a convenient way to decompose objects in
  $\Gproj\mathbb Z[\eps]$ into indecomposable objects;
  the Krull-Remak-Schmidt property, however, fails.
  The description of the indecomposable objects in
  Proposition~\ref{Gproj} does not generalize to arbitrary
  rings (for example to rings of algebraic integers)
  as we show in Remark~\ref{rem-number-rings}.
  Finally, we verify that the classification of all finitely
  generated $\mathbb Z[\eps]$-modules is not feasible.

\subsection{Gorenstein-projective modules and subgroup embeddings}
\label{sec-three-two}
We compare the category of Gorenstein-projective $\mathbb{Z}[\eps]$-modules
with the monomorphism category $\mathcal S(\proj(\mathbb{Z}))$.
This category $\mathcal S(\proj(\mathbb{Z}))$
has as objects the embeddings of a subgroup
in a finitely generated free abelian group, and morphisms are given by
commutative diagrams.
We use notion of epivalences following \cite[Definition 5.5]{GKKP2}.  
 
\begin{defin}\label{def:epiv}
  An additive functor $F:\mathcal{C}\to \mathcal{D}$ is said to be an {\bf epivalence}, if it is full, dense and reflects isomorphisms.
\end{defin} 
\medskip 
Notice that the notion of epivalences was introduced in \cite[Chapter II]{Aus} by the name of ''representation equivalence''. We prefer to use the name ``epivalence'', because it emphasizes the fact that for an epivalence $F:\mathcal{C}\to \mathcal{D}$, $\mathcal{C}$ may contain more morphisms.

\medskip 
\begin{prop}
  \label{prop-epivalence}
  There is an epivalence between $\Gproj\mathbb{Z}[\eps]$ and
  the monomorphism category $\mathcal S(\proj(\mathbb{Z}))$.     
\end{prop}

\begin{proof}
 We have seen in Lemma~\ref{isom} that each object
     $(N,e)\in \Gproj\mathbb{Z}[\eps]$ has the form $(I\oplus \Ker e,
    \lmatrix\mu)$ where
      $I$ is a subgroup of $N$ isomorphic to $\Im e$ and
      $\mu:I\to \Ker e$ is the restriction of the multiplication by $e$.
      Thus, we can assign to $(N,e)$ the inclusion
      $(\Im\mu\subset \Ker e)=(\Im e\subset\Ker e)
      \in\mathcal S(\proj(\mathbb{Z}))$.

      Moreover, a homomorphism $f:(N,e)\to (N',e')$ gives rise to
      a commutative square
      
      $$\begin{tikzcd}\Im e \arrow[r,"\rm incl"] \arrow[d,"a"] & \Ker e \arrow[d,"b"]\\
        \Im e'\arrow[r,"\rm incl"] & \Ker e',
      \end{tikzcd}$$
      hence to a morphism in $\mathcal S(\proj(\mathbb Z))$.
      This functor is full and dense, but in general,
      it is not faithful as $f$ is given by a matrix
      $(\begin{smallmatrix}a&0\\c&b\end{smallmatrix})$ where
        $c:\Ker e\to\Im e'$ can be any map,
        see Remark~\ref{rem-homomorphisms}.
  
  It is easy to check that the morphism
  $\left(\begin{smallmatrix} a&0\\c&b \end{smallmatrix}\right)$
  is an isomorphism if and only if so are $a$ and $b$.
  So $F$ reflects isomorphisms.
\end{proof}

\begin{rem} 
  An important fact about epivalences is that $F$ preserves and
  reflects indecomposable objects.
  Hence the epivalence  $F:\Gproj\mathbb{Z}[\eps]\to \mathcal S(\proj(\mathbb Z))$
  immediately yields a bijection between indecomposable objects
  in $\Gproj\mathbb{Z}[\eps]$ and $\mathcal S(\proj(\mathbb Z))$.
  That is to say, the category $\Gproj\mathbb{Z}[\eps]$
  has the same representation type as $\mathcal S(\proj(\mathbb Z))$.
\end{rem}

\subsection{Failure of Krull-Remak-Schmidt}
\label{sec-three-one}

Every Gorenstein-projective $\mathbb Z[\eps]$-module
has a decomposition as a direct
sum of indecomposable objects.  This decomposition need not be unique.
We refer to \cite{Fac} for further examples and discussion.

\begin{prop}\label{Gproj}
  Let $\mathbb Z[\eps]$ be the ring of dual integers.
  \begin{enumerate}
  \item Each Gorenstein-projective $\mathbb Z[\eps]$-module is a direct sum of
    indecomposable $\mathbb Z[\eps]$-modules.
  \item Each indecomposable Gorenstein-projective $\mathbb Z[\eps]$-module
    is isomorphic to an embedding $(I\subset \mathbb Z)$ for some ideal $I$
    in $\mathbb Z$.
  \end{enumerate}
\end{prop}

\begin{proof}
The object $(A\subset B)\in \sub(\mathbb Z[\eps])$ is given as an inclusion of a subgroup
in a finite rank free abelian group $B$, hence gives rise to a short exact
sequence
$$\mathcal E:\quad
0\longrightarrow A\longrightarrow B\stackrel g\longrightarrow C\longrightarrow 0$$
which is a projective resolution of the factor 
$C$ which has the form
$C=P\oplus T$ where $P$ is a finitely generated free $\mathbb Z$-module and $T$ is finite.
The composition $h:B\to T$ of $g$ with the canonical map $C\to P$ is a split epimorphism
  which yields a decomposition $(A\subset B)=(0\to P) \oplus (A\subset B')$
  where $B'=\Ker h$.

  \smallskip
  Clearly, the direct summand $(0\to P)$ satisfies the claims in the Proposition;
  and so does the complement $(A\subset B')$ according to the following lemma.
\end{proof}

\medskip
\begin{lem}
\label{lem-Gproj}
  Suppose $\mathcal E:0\to A\to B\stackrel h\to T\to 0$ is a short exact sequence of abelian groups
  where $A$ and $B$ are free of finite rank and $T$ is finite.  Then $\mathcal E$ is isomorphic
  to a direct sum of sequences of the form $0\to I\to \mathbb Z\to \mathbb Z/I\to 0$ where
  $I\subset \mathbb Z$ is a non-zero ideal.
\end{lem}

\begin{proof}
  Let $m\in\mathbb Z$ be such that $mT=0$, so $A$ can be considered a subgroup of $B$
  containing $mB$.  Such subgroups are in 1-1-correspondence with the subgroups of $B/mB$.
  It suffices to show that every sequence $\bar{\mathcal E}:0\to A/mB\to B/mB\to T\to 0$ is isomorphic
  to a direct sum of sequences of the form $0\to (u)/(m)\to \mathbb Z/(m)\to \mathbb Z/(u)\to 0$
  for suitable divisors $u$ of $m$. We may assume that $m=p^n$ is a prime power,
  so $B/mB\cong (\mathbb Z/(p^n))^r$ is a projective module over $\mathbb Z/(p^n)$.

  \smallskip
  
  In \cite[Theorem~5.1 and Chapter 6]{rs-art} we consider the category $\mathcal S(p^n)$ with objects the
  embeddings of a subgroup in a finite abelian $p^n$-bounded group and compute the Auslander-Reiten translates
  for the indecomposable non-projective objects; the translation is given by rotating the corresponding triangle
  in the stable category to the left.  In the special case where the embedding is $X=(U\subset P)$
  where $P$ is a projective $\mathbb Z/(p^n)$-module, the first four translates are $\tau_{\mathcal S}X=(0\subset P/U)$,
  $\tau^2_{\mathcal S}X=(P/U=P/U)$, $\tau^3_{\mathcal S}X=(P/U\subset Q)$ where $Q$ is a projective module,
  and $\tau^4_{\mathcal S}X=(0\subset U)$ (the period of $\tau_{\mathcal S}$ is six).
  Clearly, the object $(0\subset U)$ is indecomposable only if $U$ is indecomposable as a $\mathbb Z /(p^n)$-module.
  Thus, the only indecomposable objects of the form $(U\subset P)$ are the ones where $P$ itself is
  indecomposable as a $\mathbb Z/(p^n)$-module.

  \smallskip
  In particular, the sequence $\overline{\mathcal E}: 0\to A/mB\to B/mB\to T\to 0$ decomposes as indicated.
  \end{proof}
 
\begin{lem}
  The category $\Gproj\mathbb Z[\eps]$ fails to have the Krull-Remak-Schmidt
  property.
\end{lem}

\begin{proof}
  As example, we show that
$(6\mathbb Z\subset \mathbb Z)\oplus (\mathbb Z\subset\mathbb Z)$
and $(2\mathbb Z\subset \mathbb Z)\oplus (3\mathbb Z\subset \mathbb Z)$
are indecomposable decompositions of isomorphic modules.

$$
  \begin{tikzcd}[ampersand replacement=\&]
    0 \arrow[r] \& 6\mathbb Z\oplus \mathbb Z \arrow[r, "incl"]
    \& \mathbb Z\oplus \mathbb Z \arrow[r]
    \arrow[d,bend left=20, "{\begin{smallpmatrix} 1 & 2 \\ 1 & 3 \end{smallpmatrix}}"]
    \& \mathbb Z_6\oplus 0\arrow[r] \& 0 \\
    0\arrow[r] \& 2\mathbb Z\oplus 3\mathbb Z\arrow[r,"incl"]
    \& \mathbb Z\oplus \mathbb Z \arrow[r]
    \arrow[u,bend left=20, "{\begin{smallpmatrix} 3 & -2 \\ -1 & 1 \end{smallpmatrix}}"]
    \& \mathbb Z_2\oplus \mathbb Z_3\arrow [r] \& 0
  \end{tikzcd}
$$
The vertical maps preserve the subspaces and hence induce an isomorphism
  between complexes.
  However, the object $(6\mathbb Z\subset \mathbb Z)$
  is not isomorphic to either $(2\mathbb Z\subset \mathbb Z)$ or
  $(3\mathbb Z\subset \mathbb Z)$ (since the cokernels are not isomorphic as
  $\mathbb Z$-modules).
\end{proof}

\begin{rem}\label{multiplicity}
  The indecomposable summands in a direct sum decomposition of $(A\subset B)$
  are not determined uniquely, as shown above.
  Only the number of indecomposable direct summands of
  $(A\subset B)$ is given as $\rank B=\dim B\otimes \mathbb Q$ and the multiplicity of $(0\to \mathbb Z)$
  as a direct summand of $(A\subset B)$ is determined as $m=\dim B/A\otimes \mathbb Q$.
\end{rem}

  \begin{rem}
    \label{rem-number-rings}
    Many of the results in this paper, in particular all statements in
    Sections \ref{sec-2} and \ref{sec-stable}, hold true if the ring $\mathbb Z$
    of integers is replaced by any other hereditary Noetherian integral domain $R$,
      as mentioned in the introduction.

    \medskip
    However, Proposition~\ref{Gproj} is not valid in the general case
    as the proof of Lemma~\ref{lem-Gproj} cannot be adapted.
    Here is an example.

    \medskip
    Let $K=\mathbb Q(\sqrt{-5})$ be the algebraic number field and
    $R=\mathbb Z[\sqrt{-5}]$ its ring of integers.  This
    Dedekind domain is well-known
    from algebraic number theory as unique factorization fails:
    $$6=2\cdot 3=(1+\sqrt{-5})\cdot(1-\sqrt{-5})$$
    Consider the ideal $I=(2,1+\sqrt{-5})$ in $R$ generated by $2$
    and $1+\sqrt{-5}$. The norm of $I$ is 
    $$N(I)=|R/I|={\rm det}\smallmat2011=2.$$
    Since the norm of an element $x=a+b\sqrt{-5}$ is $N(x)=a^2+5b^2$,
    no element can have norm $2$, hence $I$ is not a principal ideal.
    As a consequence, the following two embeddings are not equivalent:
    $$((2)\subset I) \quad\text{and}\quad (I\subset R).$$
    Thus, the corresponding $R[\eps]$-modules cannot be isomorphic.
    But they are quasi-isomorphic as we are going to show.
    Since $2$ ramifies as $(2)=I^2$, both embeddings
    have quotient isomorphic to $R/I\cong \mathbb Z/(2)$.
    This phenomenon may be strange when dealing with finite length modules,
    hence we give an illustration.  There are $f$, $g$ which make the diagram commutative.
    $$
    \begin{tikzcd}[ampersand replacement=\&]
      0 \arrow[r] \& (2) \arrow[r, "incl"] \arrow[d,bend left=20]
      \& (2,1+\sqrt{-5}) \arrow[r, "can"]
      \arrow[d,bend left=20, "g"]
      \& \mathbb Z/(2)\arrow[r]\arrow[d,equal] \& 0 \\
      0\arrow[r] \& (2,1+\sqrt{-5}) \arrow[r,"incl"]\arrow[u,bend left=20]
      \& \mathbb Z[\sqrt{-5}] \arrow[r,"can"]
      \arrow[u,bend left=20, "f"]
      \& \mathbb Z/(2)\arrow [r] \& 0
    \end{tikzcd}
    $$
    Namely, for $x\in\mathbb Z[\sqrt{-5}]$, define $f(x)=(1+\sqrt{-5})\cdot x$
    and for $y=a+b\sqrt{-5}\in\mathbb Z[\sqrt{-5}]$ with $a\equiv b\mod 2$,
    put $g(y)=\frac12(1-\sqrt{-5})\cdot y$. Note that $(f\circ g)(y)=3y$
    and $(g\circ f)(x)=3x$.  In this sense, $f$ and $g$ are
    quasi-isomorphisms.

    \medskip
    In summary, for $R=\mathbb Z[\sqrt{-5}]$, 
    the object  $X=(N,e)=((2)\oplus I,\smallmat 00{\rm incl}0)$
    is an indecomposable
    Gorenstein-projective $R[\eps]$-module by Proposition~\ref{prop-subR}.
    Since $\Ker e=I\not\cong R$, the object $X$ does not satisfy
    condition (2) in Proposition~\ref{Gproj}.
  \end{rem}

\subsection{The category of finitely   generated $\mathbb Z[\eps]$-modules is wild}
\label{sec-three-three}

In Proposition~\ref{prop-subR} we gave a classification of all Gorenstein-projective
$\mathbb Z[\eps]$-modules. Is it possible to classify all finitely generated
$\mathbb Z[\eps]$-modules?  Or at least the finite $\mathbb Z[\eps]$-modules?
That is, the embeddings $(A\subset B)$ where $B$ is a finite abelian group
and $A$ a subgroup?

\medskip
The last problem is the well-known Birkhoff Problem; see \cite{k} for a survey.
In particular, the classification of all embeddings $(A\subset B)$ where $B$ is
a finite $p^7$-bounded group is not considered to be feasible:

\begin{rem}
  \label{rem-wild}
  For any prime number $p$ and any integer $n\geq7$,
  the category $\mathcal S(p^n)$ of all embeddings $(A\subset B)$ of a subgroup $A$ in a
  finite abelian $p^n$-bounded group $B$ is controlled wild \cite{rs-wild}.
\end{rem}
\vskip10pt
\section{The stable category of Gorenstein-projective modules}
\label{sec-stable}
\vskip10pt
In this section, we will prove Theorem~\ref{thm-main}
although some of our results are stated for a more general setting.
\vskip10pt
Let $R$ be an arbitrary hereditary Noetherian integral domain 
and $D^b(R)$ be the bounded derived category.
This is a triangulated category and its shift functor will be denoted by $[1]$.
One of our main results is that  there is an equivalence
between the orbit category $D^b( R )/[1]$ and the stable category
$ \usub(R[\eps])$.
We will prove it in Section~\ref{sec-four-two}.

\subsection{ The stable category}
\label{sec-four-one}
 
First of all, we will discuss some properties of the stable category of
$\usub(R[\eps])$.

\medskip

\begin{lem}\label{homotopy-zero}
Let $R$ be a hereditary Noetherian integral domain 
and $\proj(R)$ the subcategory consisting of the projective objects in $\mod R$.
  Let $f$ be a morphism in $ { C ^b(\proj(R))/[1]}$.
  If $f$ is homotopic to zero, then $\eta_1 (f)$, defined in Section~\ref{sec-two-four},
  factors through a projective module.
\end{lem}

\begin{proof}
  From the proof of Lemma \ref{orbit-diff}, a morphism between the indecomposables
  $X^{\bullet}$ and $Y^{\bullet}$ in 
  $C^b(\proj(R))/[1]$  has the form  $f=(a, b)\oplus (c,0)$,
  where $a: X^1\to Y^1, b:X^2\to Y^2$ and $c:X^1\to Y^2$.
  If it is homotopic to zero, then there exist 
  $$t: \ X^2\to Y^1,\ e:\ X^1\to Y^1 \;\text{and}\; h:\ X^2\to Y^2$$
  such that $a=td^1_X, b=d^1_Y t$ and $c=hd^1_X+d^1_Ye.$ 
  Since 
  $\eta_1 (X^{\bullet})=X^1\oplus X^2, \eta_1(Y^{\bullet})=Y^1\oplus Y^2$,
  we can check that $\eta_1(f)$ factors through the projective object
  $$(X^1 \oplus X^2\oplus Y^1\oplus Y^2, \smallmat 00{\smallmat et0h}0)$$
  via the morphisms
  $$\left(\begin{smallmatrix}\smallmat 00{d_X^1}0 \\
    \smallmat et0h \end{smallmatrix}\right):
  (X^1 \oplus X^2, \smallmat 00{d_X^1}0 )\to
  (X^1 \oplus X^2\oplus Y^1\oplus Y^2,
  \left(\begin{smallmatrix}0&0\\ \smallmat et0h &0\end{smallmatrix}\right))$$
    and
    $$\left(\begin{smallmatrix}\smallmat et0h,&\smallmat 00{d_Y^1}0
    \end{smallmatrix}\right):
    (X^1 \oplus X^2\oplus Y^1\oplus Y^2,
    \left(\begin{smallmatrix}0&0\\  \smallmat et0h&0\end{smallmatrix}\right))
      \to (Y^1 \oplus Y^2, \smallmat 00{d_Y^1}0).$$
\end{proof}

\medskip

By the above Lemma \ref{homotopy-zero}, we know that the full and dense functor
$$\eta_1:\ C^b(\proj(R))/[1] \to  \sub(R[\eps])$$ induces a full and dense
functor to the stable category
  $$\eta_{\mathcal P, 1}:\ K^b(\proj(R))/[1] \to \usub(R[\eps]).$$

\begin{lem}\label{K-diff-equ} The induced functor
  $ \eta_{\mathcal P, 1}:\ K^b(\proj(R))/[1] \to  
  \usub(R[\eps]) $ is faithful,
  hence an equivalence of categories.
\end{lem}

\begin{proof}
  Let
$$X^{\bullet}:\ \cdots \to 0\to X^1\stackrel{d^{1}_X }\longrightarrow X^2\to 0\to \cdots ,$$
$$ Y^{\bullet}: \ \cdots \to 0\to Y^1\stackrel{d^{1}_Y }\longrightarrow Y^2\to 0\to \cdots$$ be two indecomposables in $K^b(\proj(R))$, and $$(a, b)\oplus (c,0)\in {\rm{Hom}}_{K^b(\proj(R))/[1]}(X^{\bullet}, Y^{\bullet} ).$$   If
  $ \eta_{\mathcal P, 1}((a, b)\oplus (c,0))=\left( \begin{smallmatrix} a\ &0 \\ c\ &b\end{smallmatrix}\right)$
    factors through a projective object, say
    $(P\oplus P,\smallmat 0010)$, then we have a commutative diagram:
$$
\begin{tikzcd}
  X^1\oplus X^2 \arrow[rr, "\smallmat 00{d_X^1}0 "]
  \arrow[dd,"\smallmat a0cb"]
  \arrow[dr,"\smallmat efgh" near end]
  && X^1\oplus X^2\arrow[dd,"\smallmat a0cb" near start]
  \arrow[dr,"\smallmat efgh" ] \\
  & P\oplus P\ar[rr, "\smallmat 0010",near start,crossing over]\ar[dl, "\smallmat xyzw" near start ]
  &&P\oplus P\arrow[dl,"\smallmat xyzw"]\\
Y^1\oplus Y^2\arrow[rr,"\smallmat00{d_Y^1}0" '] && Y^1\oplus Y^2
\end{tikzcd}
$$
and $a=xhd^1_X, \ b=d^1_Y xh$ and $c=xg+zh$. Hence $(a, b)$ and $(c,0)$ are homotopic to zero,
and so is $(a, b)\oplus (c,0)$.
\end{proof}

\subsection{ The equivalence between the stable and orbit categories}
\label{sec-four-two}
 Notice the fact that the bounded derived category $D^b(\mod R)$ is equivalent to the bounded homotopy category
$K^b(\proj(R))$ for an arbitrary category $R$ with finite global dimension.
We deduce from Lemma~\ref{K-diff-equ} 
the following result:
 \begin{thm}\label{derived-diff}
   For an arbitrary hereditary Noetherian integral domain 
     $R$, there is an equivalence
  between the orbit category $D^b( \mod R)/[1]$
  and the stable category $ \usub(R[\eps])$.
  \end{thm}
 
\medskip
As a consequence of Theorem~\ref{four-cats-equ}
and Theorem~\ref{derived-diff} we obtain
$$\uGproj\mathbb{Z}[\eps] \;\cong\;\usub(\mathbb{Z}[\eps])
\;\cong \;\mathcal D^b(\mathbb{Z} )/[1],$$
which completes the proof of Theorem~\ref{thm-main}.

\section{The cube}
\label{sec-cube}

We return to the cube presented in the introduction.  To make the construction
more transparent, we attach to each vertex two categories.  One, which indicates
how the category is obtained from the category at the bottom (which is $C^b$); the
other category is the equivalent category listed in the version of the cube
in the introduction.

$$
\begin{tikzcd}
  & H^\bullet C^b/(\mathcal P,[1])\cong  \Ab_{\fin} & \\  & & \\
  C^b/(\mathcal P,[1])\cong \mathcal D^b(\mathbb Z)/[1] \arrow[uur,"H^{\bullet}_{ \mathcal P,1}"]
  & H^\bullet C^b/\mathcal P\cong \amalg \Ab_{\fin} \arrow[uu,"S_{H, \mathcal P}",near start]
  \arrow[from=ddr,"P_H"',near start]
  & H^\bullet C^b/[1]\cong \mod\mathbb Z \arrow[uul,"P_{H,1}"'] \\
  & & \\
  C^b/\mathcal P\cong K^b(\mathbb Z) \arrow[uu,"S_{\mathcal P}"]
  \arrow[uur,"H^\bullet_{\mathcal P}", near start]
  & C^b/[1]\cong \Gproj\mathbb Z[\eps] \arrow[uul,"P_1"',near start,crossing over]
  \arrow[uur,"H^\bullet_1",near start,crossing over]
  & H^\bullet C^b\cong \amalg \mod \mathbb{Z} 
  \arrow[uu,"S_H"'] \\ & &\\
  & C^b=C^b(\proj(\mathbb Z)) \arrow[uul,"P"] \arrow[uu,"S"]
  \arrow[uur,"H^{\bullet}"'] &
\end{tikzcd}
$$

We have already considered most of the categories in the cube.
It remains to mention that the stable category $C^b/\mathcal P$
is equivalent to the homotopy category $K^b(\mathbb Z)$ (see \cite{N}). 
We review some basic properties of the functors involved and
show that all six faces of the cube are commutative.

\subsection{The setting}
\label{sec-five-one}
 
Let $\proj(\mathbb Z)$ be the projective subcategory of finitely generated $\mathbb Z$-modules, and
$C^b=C^b(\proj(\mathbb Z))$ the bounded complexes over $\proj(\mathbb Z)$.  Morphisms in $C^b$ are
commutative diagrams. 
By $[1]:C^b\to C^b$ we denote the usual shift (to the left).

\medskip
Consider the special object
$$P^\bullet:\cdots\to0\to \mathbb Z=\mathbb Z\to 0\to \cdots$$ (concentrated in degrees $1,0$)
and its shifts are the indecomposable projective objects in $C^b$; we denote by
$\mathcal P$ the categorical ideal in $C^b$ of all maps which factor through a
projective object in $C^b$.

  \medskip
  Since $P^\bullet$ and its shifts generate the subcategory of all projective objects in $C^b$,
  by Theorem \ref{four-cats-equ} we have
  
  \begin{lem} There is an equivalence between $C^b/(\mathcal P,[1])$
    and $\uGproj\mathbb{Z}[\eps]$.
  \end{lem}

  \subsection{The functors in the cube}
  \label{sec-five-two}
First, we consider the push-down  functors
$$S:C^b\to C^b/[1],\quad\text{and}\quad
S_{\mathcal P}: C^b/\mathcal P\to C^b/(\mathcal P,[1]).$$

\begin{lem}
  The push-down functors are dense and faithful.
\end{lem}

\begin{proof}
  It is clear that the push-down functors are faithful.
  For the density, it suffices to consider the functor $S$.
  We have seen in Proposition \ref{prop-cb} that each indecomposable object in the target
  category occurs as the image of an object of width at most two.
\end{proof}

There are two more push-down functors,
$$S_H:\amalg \mod \mathbb Z\to\mod\mathbb Z\quad\text{and}\quad
S_{H,\mathcal P}:\amalg \Ab_{\fin}\to \Ab_{\fin},$$
which are given by taking the direct sum of the objects and morphisms.

\medskip
In the cube, the four push-down functors are represented by arrows pointing up.

\bigskip

Next, we consider the four
canonical functors modulo the
(categorical) ideal $\mathcal P$ of all maps which factor through a
projective object:
$$P:C^b\to C^b/\mathcal P, \quad
P_1:C^b/[1]\to C^b/(\mathcal P,[1]),\quad
P_{H}: H^\bullet C^b\to H^\bullet C^b/\mathcal P$$ and
$P_{H,1}: H^\bullet C^b/[1]\to H^\bullet C^b/(\mathcal P,[1])$,
are represented in the cube by arrows towards the left.
Obviously,

\begin{lem} Each canonical functor modulo projectives is full and dense.\qed
\end{lem}

\bigskip
 
The remaining four functors are given by taking homology.
They are indicated in the cube by arrows pointing towards the right.  The functor
$$ H^{\bullet}:C^b \to H^\bullet C^b$$
sends each bounded complex
$\cdots\to X^i\stackrel{d^i} \to X^{i+1}\stackrel{d^{i+1}}\to X^{i+2}\stackrel{d^{i+2}}\to \cdots$
to the object $ ( \ker d^{i+1}/\Ima d^i)_{i\in \mathbb Z}\in\amalg\mod\mathbb Z$.
The action of $H^{\bullet}$ on morphisms is induced by the property of cohomology.
Note that $H^\bullet$ maps projective objects to zero, hence factors through $P$:
$$H^\bullet=\uH^\bullet\circ P.$$
Hence we can define $H^\bullet_{\mathcal P}=P_H\circ \uH^\bullet$.

\smallskip
The functor
$H^{\bullet}_1: C^b/[1]\to H^\bullet C^b/[1]$
sends an indecomposable object 
$\cdots\to 0\to U\stackrel{u}\hookrightarrow F\to 0\to \cdots $
to its homology $F/U\in\mod\mathbb Z$.
As $H^\bullet$, also $H^\bullet_1$ maps projective objects to zero, hence factors through $P_1$,
$$H^\bullet_1=\uH^\bullet_1\circ P_1,$$
and we can define $H^\bullet_{\mathcal P,1}=P_{H,1}\circ\uH^\bullet_1$.
We have as before:

\begin{lem}
  The homology functors above are full and dense.\qed
\end{lem}

\subsection{The faces of the cube are commutative}
\label{sec-five-three}
There are two squares given by the functors $S$ and $P$ which
form the front left and rear right faces 
of the cube:
$$
\begin{tikzcd}  C^b \arrow[r,"P"] \arrow[d,"S"] & C^b/\mathcal P \arrow[d,"S_{\mathcal P}"]\\
  C^b/[1] \arrow[r,"P_1"] & C^b/(\mathcal P,[1])
\end{tikzcd}
\qquad
\begin{tikzcd}
 \amalg \mod \mathbb Z\arrow[r,"P_H"] \arrow[d,"S_H"]
  & \amalg \mathcal{A}b_{\rm fin}  \arrow[d,"S_{H,\mathcal P}"]\\
 \mod \mathbb Z \arrow[r,"P_{H,1}"] & \mathcal{A}b_{\rm fin}
\end{tikzcd}
$$
Both squares are commutative since $S$ and $S_H$ preserve  projectives.

\medskip
 
Next, consider the squares given by the functors $H^\bullet,P$,
which form the top and bottom faces of the cube:
$$
\begin{tikzcd}
  C^b \arrow[r,"P"] \arrow[d,"H^{\bullet}"] & C^b/\mathcal P \arrow[d,"H^{\bullet}_{\mathcal P}"]\\
  \amalg \mod \mathbb Z
  \arrow[r,"P_H"] & \amalg \Ab_{\fin}
\end{tikzcd}
\qquad
\begin{tikzcd}
  C^b/[1] \arrow[r,"P_1"] \arrow[d,"H^{\bullet}_1"] & C^b/(\mathcal P,[1]) \arrow[d,"H^{\bullet}_{\mathcal P,1}"]\\
  \mod\mathbb Z \arrow[r,"P_{H,1}"] &  \Ab_{\fin}.
\end{tikzcd}
$$
Both squares are commutative by definition of the functors $H^\bullet_{\mathcal P}$
and $H^\bullet_{\mathcal P,1}$:
$$H^\bullet_{\mathcal P}\circ P=P_H\circ\uH^\bullet\circ P=P_H\circ H^\bullet,\;\text{and}\;
H^\bullet_{\mathcal P,1}\circ P_1=P_{H,1}\circ\uH^\bullet_1\circ P_1=P_{H,1}\circ H^\bullet_1$$

\medskip
Finally, the squares given by the functors $S, H^\bullet$
form the front right and rear left squares of the cube:

$$\begin{tikzcd}  C^b \arrow[r,"S"] \arrow[d,"H^{\bullet}"] & C^b/[1] \arrow[d,"H^{\bullet}_1"]\\
  \amalg \mod \mathbb Z
\arrow[r,"S_H"] &  \mod \mathbb Z
\end{tikzcd}
\qquad
\begin{tikzcd}  C^b/\mathcal P \arrow[r,"S_{\mathcal P}"] \arrow[d,"H^{\bullet}_{\mathcal P}"] &
  C^b/({\mathcal P,[1]})\arrow[d,"H^{\bullet}_{\mathcal P,1}"]\\
  \amalg \Ab_{\fin}\arrow[r,"S_{\mathcal P,1}"] & \Ab_{\fin}
\end{tikzcd}$$
The squares are commutative since cohomology commutes with the push-down functors.

\medskip
In conclusion, we have:

\begin{prop}
  All faces in the cube are commutative squares. \qed
\end{prop} 

  \medskip
  We hope that the categories and functors in the cube help to make the
  interplay between the categories of bounded complexes $C^b(\proj(\mathbb Z))$,
  of Gorenstein-projective modules over the dual integers $\Gproj\mathbb Z[\eps]$,
  of orbits in the derived category of the integers $\mathcal D^b(\mathbb Z)/[1]$,
  and of finite abelian groups $\Ab_{\fin}$ a bit more transparent.

\bigskip
  {\bf Acknowledgement.}  The authors would like to thank Shijie Zhu for advice which has led
  to  substantial improvements of the paper.
 They are grateful to Jos\'e V\'elez-Marulanda for pointing out
    the application to rings of algebraic integers.

\end{document}